\documentclass[iicol,sn-mathphys,Numbered]{sn-jnl}

\usepackage{graphicx}%
\usepackage{multirow}%
\usepackage{amsmath,amssymb,amsfonts}%
\usepackage{amsthm}%
\usepackage{mathrsfs}%
\usepackage[title]{appendix}%
\usepackage{xcolor}%
\usepackage{textcomp}%
\usepackage{manyfoot}%
\usepackage{booktabs}%
\usepackage{algorithm}%
\usepackage{algorithmicx}%
\usepackage{algpseudocode}%
\usepackage{listings}%

\usepackage{float}
\usepackage{amssymb}
\usepackage{tikz}
\usepackage{amsthm,amsfonts,amsmath}
\usepackage[utf8]{inputenc}
\usepackage{cleveref}
\usepackage{todonotes}
\usepackage{verbatim}

\theoremstyle{thmstyleone}%
\newtheorem{theorem}{Theorem}
\theoremstyle{thmstyletwo}%

\theoremstyle{thmstylethree}%

\newcommand{\T}{\mathcal{T}}

\newcommand{\csum}{\sum_{j} C_j}

\newcommand{\R}{\mathbb{R}}
\newcommand{\tequal}{t^<}

\newcommand{\lefttodo}{\overline{W}}
\newtheorem{claim}{Claim}
\newtheorem{lemma}{Lemma}

\newtheorem{corollary}{Corollary}
\newtheorem{observation}{Observation}
\newcommand{\ignore}[1]{}
\newenvironment{claimproof}[1]{\par\noindent\emph{Proof.}\space#1}{\hfill~$\blacksquare$\\}

\begin{document}

\title[Scheduling jobs that change over time]{Scheduling jobs that change over time}


\author[2]{\fnm{Roel} \sur{Lambers}}\email{r.lambers2@hva.nl}
\equalcont{These authors contributed equally to this work.}

\author[1]{\fnm{Rudi} \sur{Pendavingh}}\email{r.a.pendavingh@tue.nl}
\equalcont{These authors contributed equally to this work.}

\author[1]{\fnm{Frits} \sur{Spieksma}}\email{f.c.r.spieksma@tue.nl}
\equalcont{These authors contributed equally to this work.}

\author*[1]{\fnm{C\'eline M.F.} \sur{Swennenhuis}}\email{cmfswennenhuis@gmail.com}
\equalcont{These authors contributed equally to this work.}

\affil[1]{\orgdiv{Department of Mathematics and Computer Science}, \orgname{Eindhoven University of Technology}, \orgaddress{\street{PO Box 513}, \city{Eindhoven}, \postcode{5600MB}, 
\country{the Netherlands}}}

\affil[2]{\orgname{Amsterdam University of Applied Sciences}, 
\city{Amsterdam}, 
\country{the Netherlands}}




\abstract{We consider a 1-machine scheduling problem where the temperature of a job rises during processing, and cools down when not being processed according to given linear heating and cooling rates. No job's temperature is allowed to rise above a given threshold, and no job's temperature can cool below 0. Another crucial property of our problem is that jobs can be preempted an arbitrary number of times, and even more, we allow that a job is processed for an infinitely small amount of time. We consider two objectives: minimize the makespan, and minimize the sum of completion times. Our results are as follows. We show how to compactly represent a solution. Further, we prove that the problem of minimizing the sum of completion times can be solved in polynomial time by formulating it as a linear program, and deriving a structural property. This result can be extended to hold for any number of machines. Further, we show that a minimum makespan can be found in~$O(n)$ time, even when heating and cooling rates are job-dependent.}




\maketitle

\section{Introduction}
\label{sec:intro}
Traditionally, scheduling problems are formulated in terms of machines and jobs. The machines represent resources, and the jobs represent tasks to be carried out by the machines. In the vast majority of cases, the jobs are characterized by a set of properties (processing time, due dates/deadlines, weights, ...) that do not change while the job is being processed. In this paper however, we investigate a setting where a relevant property of a job changes over time as a function of being processed, and not being processed. We denote this property by ``temperature''; notice however that this choice of words stems from a particular application, see Section~\ref{sec:motivation} for further details.

More specifically, in our setting, the temperature of a job rises during processing, and cools down when not being processed. In addition, no job's temperature should rise above a given threshold, and no job's temperature can cool below 0.  
Another crucial property of our problem is that we allow for an indefinite number of preemptions of a job, and even more, we allow that a job is processed for an arbitrarily small amount of time --- we defer a precise problem statement to Section~\ref{sec:probstat}.
These properties give rise to a class of scheduling problems that  has hitherto received little to no attention. 

We consider two traditional scheduling objectives: the sum of completion times, and the makespan. We show that each of these problems can be solved in polynomial time.

\subsection{Motivation}
\label{sec:motivation}
The following situation is described in Pham et al.~\cite{Phametal2020}. Given is a 2-dimensional piece of material with locations that need to be ``processed''. Processing a location means that a laser (or some other device) is used to heat up the given location for a given amount of time in order to make it suitable for further processing. However, it may well be the case that one cannot heat a job continuously for as long as its processing time, as the temperature on that location would rise too high (and the material would melt, and become unusable). Thus, preemptions are necessary to allow the material to cool, after which processing again becomes possible. Notice that the resulting problem is more general than the one we will study here, as Pham et al.~\cite{Phametal2020} also take travel times between locations into account.

Another situation that motivates our problem can be found in the process industry (see Section~\ref{sec:lit} for references). Imagine a number of ``liquid-producing'' devices, each equipped with a button; there is also a buffer associated with each device. The machine needs to select a button in order for the corresponding device to produce its liquid. When some device's button is selected, the liquid enters the device's buffer with a particular rate (this is akin to the heating rate). As the buffer has a particular capacity, there is no point in continuing to select this button when the buffer is full. When a machine's button is not selected, its buffer empties with a particular rate (the cooling rate). The amounts of each liquid to be produced determine the processing times associated with each device.

\subsection{Related Literature}
\label{sec:lit}

We describe a number of problem settings that share characteristics with our problem.

First, when considering the challenge of controlling temperatures in a job scheduling context, there is quite some literature on operating a micro-processor. A micro-processor carries out computational tasks, and its temperature is determined, among other factors, by the speed with which it carries out these tasks. This speed determines how fast the tasks can be completed. We mention the following papers dealing with this subject: Bampis et al.~\cite{Bampisetal2013}, Birks and Fung~\cite{BirFun2013, BirFun2017}, and Chrobak et al.~\cite{Chrobaketal2011}. A key difference in these works compared to this contribution is that the temperature of the {\em system}, i.e., the temperature of the machine is of prime concern, whereas in our case, the temperature of each job matters. Bai et al.~\cite{Baietal2008} describe a one machine scheduling problem where each job has a given processing time, and all jobs cool down during processing according to an exponential function. The objective here is to schedule the jobs so that total temperature loss is minimized. They establish the complexity of different cases of the problem. Our problem is fundamentally different from this setting, as (i) in our case the temperature of a job changes also when not processed, (ii) we allow preemptions.

Secondly, Martinelli and Valigi~\cite{Marval1998} describe a production control problem where a single machine can produce parts of different types. There is a buffer corresponding to each part type, and there is known, constant demand rate for each part type. The problem is to decide, at each moment in time, which part type to produce, and at which rate; the objective in this setting is to minimize total (discounted) inventory costs which relate to the number of parts in each buffer. The demand rate corresponds to the cooling rate, and the rate at which the machine produces parts of a certain type, corresponds to the heating rate. Martinelli and Valigi~\cite{Marval2002, Marval2004} establish, for the case of two part types, optimal policies in the presence of setup times, backlog, and capacitated buffers. There are, however, fundamental differences between our problem and this production control setting:  the objective is fundamentally different, there is no analog of backlog in our problem, and it is a fact that in our setting a job either cools down or heats up - which is not true in production control settings, as the production rate can be chosen, and is independent of the demand rate.

More generally, we can connect our problem to a class of problems known as {\em polling systems}. In a polling system a set of queues (jobs in our terminology) is served by a single server (a machine) that can serve at most one queue at a time. Typically, in a polling system, when a queue is not served, its length (temperature) increases according to a given arrival process (when this arrival process boils down to a constant deterministic arrival rate, this would correspond with the cooling rate). When the server is serving a queue, there is a known rate with which the queue empties (corresponding to the heating rate). Typical objectives focus on either the total queue length in the system, or focus on waiting times. Many variations of polling systems have been investigated, and there is an immense variety of applications. Here, we are content to refer to the following survey papers that describe the existing results and the applications of polling systems: Boon et al.~\cite{boonetal2011}, Vishnevskii and Semenova~\cite{Vissem2006}, Borst and Boxma~\cite{borbox2018}, Takagi~\cite{Tak2000}. We also mention Kruskal~\cite{kruskal1968} for an early contribution in a deterministic setting, and Lefeber et al.~\cite{leflamroo2011} and Matveev et al.~\cite{matetal2016} who also deal with a deterministic variant of a polling system. Deterministic multi server problems are addressed in Bertsimas and Kim~\cite{berkim2023}.

Our problem, while sharing similarities with polling systems, is different from problems analyzed in the literature on polling systems. Indeed, since the queue length in a polling system corresponds to the temperature in our scheduling problem, serving a queue means, in our context, actually enlarging the queue length. Other differences are that: (i) our input is completely deterministic, and (ii) our objective does not arise from the length of the queues, or the waiting times.

Finally, as described in Section~\ref{sec:motivation}, we mention that our problem can be relevant when dealing with liquids, where the temperature of a job represents the relative used capacity of a buffer or a tank. The heating (cooling) rate is then interpreted as the rate with which a buffer is filled (empties). While the application sketched in Section~\ref{sec:motivation} is highly stylized, more complicated practical situations involving the handling of liquids have been reported in literature. Indeed, in the process industry there is ample literature on the operational challenges that arise when producing liquids, see Kallrath~\cite{Kallrath2002} and Kilic~\cite{Kilic2011} for surveys in this field. One noteworthy case is the so-called pooling problem (see Audet et al.~\cite{Audetetal2004} and Gupte et al.~\cite{Gupteetal2017}); this problem is a well-studied generalization of min-cost flow. Given suppliers with raw materials containing known specifications, the problem is to find a minimum cost way of mixing these materials in intermediate tanks (pools) so as to meet the demand requirements and their specifications. While this problem focuses on mixing various inflows, our problem could arise when given very specific network characteristics.

\section{Problem Statement}
\label{sec:probstat}
Given is a single machine and a set of jobs~$J$ that need to be processed by the machine. The machine can process one job at a time. There is a given, positive processing time \(p_j\)  associated to each job~$j \in J$. Jobs can be pre-empted.

For each time~$t \geq 0$ and for each job~$j \in J$, there is a temperature \(T_j(t)\) at time \(t\), with~$T_j(0)=0$. The temperature of a job changes over time according to the following rules:
\begin{itemize}
    \item When a job~$j \in J$ is not processed, it cools with factor \(\alpha\) (notice that \(\alpha < 0\)).
    \item When a job~$j \in J$ is processed, it heats with factor \(\beta\) (notice that \(\beta > 0\)).
    \item A job~$j \in J$ can never have a temperature below \(0\). If \(T_j(t) = 0\) for some~$t \geq 0$ and some~$j \in J$, then job~$j$ will not cool down any further, even when not processed.
\end{itemize}
The temperature of a job is not allowed to exceed a certain threshold; we choose this threshold equal to \(1\) for all jobs (we can make this assumption, as we can scale~$\alpha$ to~$\alpha /T$ and~$\beta$ to~$\beta/T$).  Thus, an instance of our problem is completely specified by~$\alpha, \beta$, and~$p_j$ for~$ j \in J$; we assume that all these data are rational.

 We focus on the following two objectives: minimizing the sum of completion times (Section~\ref{sec:sumofcompletiontimes}) and minimizing the makespan (Section~\ref{sec:makespan}). It is not a priori clear that these minima are attained by schedules in the sense described above, and so our goal is to find the infimum value in each case.

We refer to the former problem, minimizing the sum of completion times,  as problem~$Q1^{sum}$, and to the latter problem, minimizing the makespan, as~$Q1^{max}$. 
We will also consider the generalization that arises when~$m$ machines are available; these problems are denoted by~$Qm^{sum}$ and~$Qm^{max}$.

Our contribution is as follows.
\begin{itemize}
\item We show how to compactly represent a solution of our problems (Section~\ref{sec:prepinsight}).
\item We prove that problem~$Q1^{sum}$ can be solved in polynomial time by (i) 
formulating a linear program that needs as input a sequence of jobs, and (ii) showing that there exists an optimum solution in which jobs are sequenced according to their processing times. This result can be extended to hold for any number of machines, i.e., for~$Qm^{sum}$ (Section~\ref{sec:sumofcompletiontimes}).
\item We show that problem~$Qm^{max}$ can be solved in~$O(n)$ time even when heating rates and cooling rates are job-dependent (Section~\ref{sec:makespan}).
\end{itemize}

\subsection{An example}\label{ssec:example}
Consider the following instance consisting of a single machine and two jobs (i.e.,~$J=\{1,2\}$) with~$p_1 = p_2 = 2, \alpha = -\frac13, \beta = 1$. What can we say about the makespan~$C_{\max}$ and sum of completion times~$\sum_{j \in J}C_j$ of this instance?

Consider what happens if we only want to process job 1. The following schedule is then optimal: processing from~$t=0$ until~$t=1$, and from~$t=4$ until~$t=5$, while not processing (i.e. cooling) in the remaining time intervals. 
Clearly, the original instance cannot be completed faster than an instance consisting solely of job 1 (or solely of job 2). As a schedule that only processes job 1 already requires 5 time-units, it becomes clear that our instance needs at least 5 time-units to complete, and that the sum of completion times is at least 5+5=10.
Indeed, any schedule that starts processing job 1 or job 2 at time~$\epsilon>0$ cannot complete that job before~$5+\epsilon$, so that the makespan of such a schedule is at least~$5+\epsilon$ and the sum of completion times is at least~$10+\epsilon$.

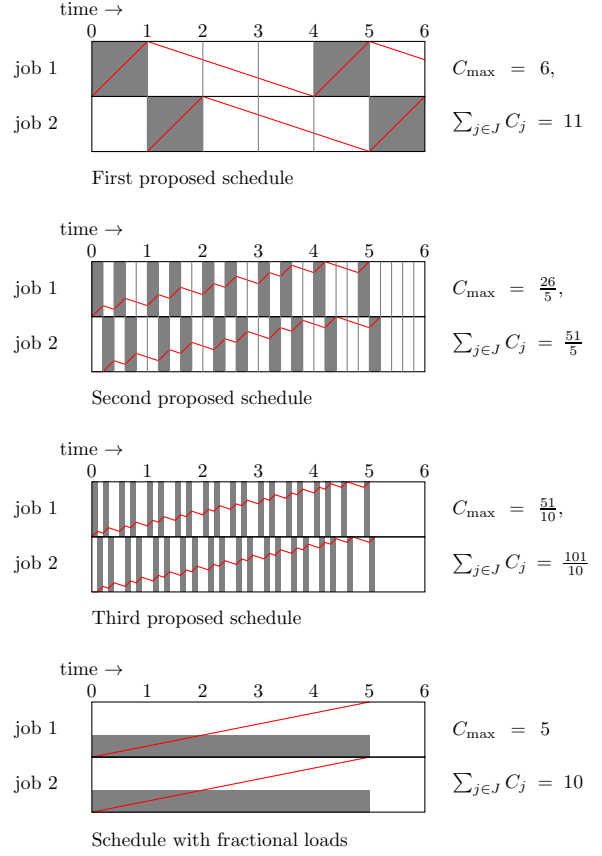
\begin{figure}[h!]
    \centering
    \begin{tikzpicture} [scale = .73, transform shape,
]

\fill[color = gray] (0,0) rectangle  (1,1);
\fill[color = gray] (4,0) rectangle  (5,1);

\draw[color = red] (0,0) -- (1,1) -- (4,0) -- (5,1) -- (6,0.667);

\fill[color = gray] (1,-1) rectangle  (2,0);
\fill[color = gray] (5,-1) rectangle  (6,0);

\draw[color = red] (0,-1) -- (1,-1) -- (2,0) -- (5,-1) -- (6,0);

\draw[color = black] (0,0) rectangle  (6,1);
\draw[color = black] (0,0) rectangle  (6,-1);
\draw[color = gray] (1,-1) -- (1,1);
\draw[color = gray] (2,-1) -- (2,1);
\draw[color = gray] (3,-1) -- (3,1);
\draw[color = gray] (4,-1) -- (4,1);
\draw[color = gray] (5,-1) -- (5,1);
 \node[] at (0,1.2) 
    {$0$};
 \node[] at (1,1.2) 
    {$1$};
 \node[] at (2,1.2) 
    {$2$};
 \node[] at (3,1.2) 
    {$3$};
 \node[] at (4,1.2) 
    {$4$};
 \node[] at (5,1.2) 
    {$5$};
 \node[] at (6,1.2) 
    {$6$};
 \node[] at (0,1.6) 
    {time $\rightarrow$};

 \node[] at (-1,.5) 
    {job $1$};
 \node[] at (-1,-.5) 
    {job $2$};

 \node[text width=6cm] at (3,-1.5) 
    {First proposed schedule};

 \node[text width=3cm] at (8,.5) 
    {$C_{\max} = 6$,};

 \node[text width=3cm] at (8,-.5) 
    {$\sum_{j \in J}C_j = 11$};

\fill[color = gray] (0,-4) rectangle  (.2,-3);
\fill[color = gray] (0.4,-4) rectangle  (.6,-3);
\fill[color = gray] (1,-4) rectangle  (1.2,-3);
\fill[color = gray] (1.4,-4) rectangle  (1.6,-3);
\fill[color = gray] (2,-4) rectangle  (2.2,-3);
\fill[color = gray] (2.4,-4) rectangle  (2.6,-3);
\fill[color = gray] (3,-4) rectangle  (3.2,-3);
\fill[color = gray] (3.4,-4) rectangle  (3.6,-3);

\fill[color = gray] (4,-4) rectangle  (4.2,-3);
\fill[color = gray] (4.8,-4) rectangle  (5,-3);

\draw[color = red] (0,-4) -- (0.2,-3.8) --	(0.4,-3.866666667) -- (0.6,-3.666666667) -- (0.8,-3.733333333) -- (1,-3.8) -- (1.2,-3.6) -- (1.4, -3.666666667) -- (1.6,-3.466666667) -- (1.8,-3.533333333) -- (2,-3.6) -- (2.2,-3.4) -- (2.4,-3.466666667) -- (2.6,-3.266666667) -- (2.8,-3.333333333) -- (3,-3.4) -- (3.2,-3.2) -- (3.4,-3.266666667) -- (3.6,-3.066666667) -- (3.8,-3.133333333) -- (4,-3.2) -- (4.2,-3) -- (4.4, -3.066666667) -- (4.6,-3.133333333) -- (4.8, -3.2) -- (5,-3)
;

\fill[color = gray] (0.4,-4) rectangle  (.2,-5);
\fill[color = gray] (0.8,-4) rectangle  (.6,-5);
\fill[color = gray] (1.4,-4) rectangle  (1.2,-5);
\fill[color = gray] (1.8,-4) rectangle  (1.6,-5);
\fill[color = gray] (2.4,-4) rectangle  (2.2,-5);
\fill[color = gray] (2.8,-4) rectangle  (2.6,-5);
\fill[color = gray] (3.4,-4) rectangle  (3.2,-5);
\fill[color = gray] (3.8,-4) rectangle  (3.6,-5);

\fill[color = gray] (4.4,-4) rectangle  (4.2,-5);
\fill[color = gray] (5.2,-4) rectangle  (5,-5);

\draw[color = red] (0,-5) -- (0.2,-5) --(0.4,-4.8) --	(0.6,-4.866666667) -- (0.8,-4.666666667) -- (1,-4.733333333) -- (1.2,-4.8) -- (1.4,-4.6) -- (1.6, -4.666666667) -- (1.8,-4.466666667) -- (2,-4.533333333) -- (2.2,-4.6) -- (2.4,-4.4) -- (2.6,-4.466666667) -- (2.8,-4.266666667) -- (3,-4.333333333) -- (3.2,-4.4) -- (3.4,-4.2) -- (3.6,-4.266666667) -- (3.8,-4.066666667) -- (4,-4.133333333) -- (4.2,-4.2) -- (4.4,-4) -- (4.6, -4.066666667) -- (4.8,-4.133333333) -- (5, -4.2) -- (5.2,-4)
;

\draw[color = black] (0,-4) rectangle  (6,-3);
\draw[color = black] (0,-4) rectangle  (6,-5);

\draw[color = gray] (0.2,-5) -- (0.2,-3);
\draw[color = gray] (0.4,-5) -- (0.4,-3);
\draw[color = gray] (0.6,-5) -- (0.6,-3);
\draw[color = gray] (0.8,-5) -- (0.8,-3);

\draw[color = gray] (1,-5) -- (1,-3);
\draw[color = gray] (1.2,-5) -- (1.2,-3);
\draw[color = gray] (1.4,-5) -- (1.4,-3);
\draw[color = gray] (1.6,-5) -- (1.6,-3);
\draw[color = gray] (1.8,-5) -- (1.8,-3);

\draw[color = gray] (2,-5) -- (2,-3);
\draw[color = gray] (2.2,-5) -- (2.2,-3);
\draw[color = gray] (2.4,-5) -- (2.4,-3);
\draw[color = gray] (2.6,-5) -- (2.6,-3);
\draw[color = gray] (2.8,-5) -- (2.8,-3);

\draw[color = gray] (3,-5) -- (3,-3);
\draw[color = gray] (3.2,-5) -- (3.2,-3);
\draw[color = gray] (3.4,-5) -- (3.4,-3);
\draw[color = gray] (3.6,-5) -- (3.6,-3);
\draw[color = gray] (3.8,-5) -- (3.8,-3);
\draw[color = gray] (4,-5) -- (4,-3);

\draw[color = gray] (4.2,-5) -- (4.2,-3);
\draw[color = gray] (4.4,-5) -- (4.4,-3);
\draw[color = gray] (4.6,-5) -- (4.6,-3);
\draw[color = gray] (4.8,-5) -- (4.8,-3);
\draw[color = gray] (5,-5) -- (5,-3);

\draw[color = gray] (5.2,-5) -- (5.2,-3);
\draw[color = gray] (5.4,-5) -- (5.4,-3);
\draw[color = gray] (5.6,-5) -- (5.6,-3);
\draw[color = gray] (5.8,-5) -- (5.8,-3);

 \node[] at (0,-2.8) 
    {$0$};
 \node[] at (1,-2.8) 
    {$1$};
 \node[] at (2,-2.8) 
    {$2$};
 \node[] at (3,-2.8) 
    {$3$};
 \node[] at (4,-2.8) 
    {$4$};
 \node[] at (5,-2.8) 
    {$5$};
 \node[] at (6,-2.8) 
    {$6$};
 \node[] at (0,-2.4) 
    {time $\rightarrow$};

 \node[] at (-1,-3.5) 
    {job $1$};
 \node[] at (-1,-4.5) 
    {job $2$};

 \node[text width=6cm] at (3,-5.5) 
    {Second proposed schedule};

 \node[text width=3cm] at (8,-3.5) 
    {$C_{\max} = \frac{26}{5}$,};

 \node[text width=3cm] at (8,-4.5) 
    {$\sum_{j \in J}C_j = \frac{51}{5}$};


\fill[color = gray] (0,-8) rectangle  (0.1,-7);
\fill[color = gray] (0.2,-8) rectangle  (0.3,-7);
\fill[color = gray] (0.5,-8) rectangle  (0.6,-7);
\fill[color = gray] (0.7,-8) rectangle  (0.8,-7);

\fill[color = gray] (1,-8) rectangle  (1.1,-7);
\fill[color = gray] (1.2,-8) rectangle  (1.3,-7);
\fill[color = gray] (1.5,-8) rectangle  (1.6,-7);
\fill[color = gray] (1.7,-8) rectangle  (1.8,-7);

\fill[color = gray] (2,-8) rectangle  (2.1,-7);
\fill[color = gray] (2.2,-8) rectangle  (2.3,-7);
\fill[color = gray] (2.5,-8) rectangle  (2.6,-7);
\fill[color = gray] (2.7,-8) rectangle  (2.8,-7);

\fill[color = gray] (3,-8) rectangle  (3.1,-7);
\fill[color = gray] (3.2,-8) rectangle  (3.3,-7);
\fill[color = gray] (3.5,-8) rectangle  (3.6,-7);
\fill[color = gray] (3.7,-8) rectangle  (3.8,-7);

\fill[color = gray] (4,-8) rectangle  (4.1,-7);
\fill[color = gray] (4.2,-8) rectangle  (4.3,-7);
\fill[color = gray] (4.5,-8) rectangle  (4.6,-7);
\fill[color = gray] (4.9,-8) rectangle  (5,-7);

\draw[color = red] (0,-8) -- (0.1, -7.9) --(0.2, -7.933333333) -- (0.3, -7.833333333) -- (0.4, -7.866666667) -- (0.5, -7.9) -- (0.6, -7.8) -- (0.7, -7.833333333) -- (0.8, -7.733333333) -- (0.9, -7.766666667)-- (1.0, -7.8)-- (1.1, -7.7)-- (1.2, -7.733333333)-- (1.3, -7.633333333)-- (1.4, -7.666666667)-- (1.5, -7.7)-- (1.6, -7.6)-- (1.7, -7.633333333)-- (1.8, -7.533333333)-- (1.9, -7.566666667)-- (2.0, -7.6)-- (2.1, -7.5)-- (2.2, -7.533333333)-- (2.3, -7.433333333)-- (2.4, -7.466666667)-- (2.5, -7.5)-- (2.6, -7.4)-- (2.7, -7.433333333)-- (2.8, -7.333333333)-- (2.9, -7.366666667)-- (3.0, -7.4)-- (3.1, -7.3)-- (3.2, -7.333333333)-- (3.3, -7.233333333)-- (3.4, -7.266666667)-- (3.5, -7.3)-- (3.6, -7.2)-- (3.7, -7.233333333)-- (3.8, -7.133333333)-- (3.9, -7.166666667)-- (4.0, -7.2)-- (4.1, -7.1)-- (4.2, -7.133333333)-- (4.3, -7.033333333)-- (4.4, -7.066666667)-- (4.5, -7.1)-- (4.6, -7)-- (4.7, -7.033333333)-- (4.8, -7.066666667) -- (4.9,-7.1) -- (5,-7)
;
	

\fill[color = gray] (0.2,-8) rectangle  (0.1,-9);
\fill[color = gray] (0.4,-8) rectangle  (0.3,-9);
\fill[color = gray] (0.7,-8) rectangle  (0.6,-9);
\fill[color = gray] (0.9,-8) rectangle  (0.8,-9);

\fill[color = gray] (1.2,-8) rectangle  (1.1,-9);
\fill[color = gray] (1.4,-8) rectangle  (1.3,-9);
\fill[color = gray] (1.7,-8) rectangle  (1.6,-9);
\fill[color = gray] (1.9,-8) rectangle  (1.8,-9);

\fill[color = gray] (2.2,-8) rectangle  (2.1,-9);
\fill[color = gray] (2.4,-8) rectangle  (2.3,-9);
\fill[color = gray] (2.7,-8) rectangle  (2.6,-9);
\fill[color = gray] (2.9,-8) rectangle  (2.8,-9);

\fill[color = gray] (3.2,-8) rectangle  (3.1,-9);
\fill[color = gray] (3.4,-8) rectangle  (3.3,-9);
\fill[color = gray] (3.7,-8) rectangle  (3.6,-9);
\fill[color = gray] (3.9,-8) rectangle  (3.8,-9);

\fill[color = gray] (4.2,-8) rectangle  (4.1,-9);
\fill[color = gray] (4.4,-8) rectangle  (4.3,-9);
\fill[color = gray] (4.7,-8) rectangle  (4.6,-9);
\fill[color = gray] (5,-8) rectangle  (5.1,-9);

\draw[color = red] (0,-9) -- (0.1,-9) -- (0.2, -8.9) --(0.3, -8.933333333) -- (0.4, -8.833333333) -- (0.5, -8.866666667) -- (0.6, -8.9) -- (0.7, -8.8) -- (0.8, -8.833333333) -- (0.9, -8.733333333) -- (1, -8.766666667)-- (1.1, -8.8)-- (1.2, -8.7)-- (1.3, -8.733333333)-- (1.4, -8.633333333)-- (1.5, -8.666666667)-- (1.6, -8.7)-- (1.7, -8.6)-- (1.8, -8.633333333)-- (1.9, -8.533333333)-- (2, -8.566666667)-- (2.1, -8.6)-- (2.2, -8.5)-- (2.3, -8.533333333)-- (2.4, -8.433333333)-- (2.5, -8.466666667)-- (2.6, -8.5)-- (2.7, -8.4)-- (2.8, -8.433333333)-- (2.9, -8.333333333)-- (3, -8.366666667)-- (3.1, -8.4)-- (3.2, -8.3)-- (3.3, -8.333333333)-- (3.4, -8.233333333)-- (3.5, -8.266666667)-- (3.6, -8.3)-- (3.7, -8.2)-- (3.8, -8.233333333)-- (3.9, -8.133333333)-- (4, -8.166666667)-- (4.1, -8.2)-- (4.2, -8.1)-- (4.3, -8.133333333)-- (4.4, -8.033333333)-- (4.5, -8.066666667)-- (4.6, -8.1)-- (4.7, -8)-- (4.8, -8.033333333)-- (4.9, -8.066666667) -- (5,-8.1) -- (5.1,-8)
;

\draw[color = black] (0,-8) rectangle  (6,-7);
\draw[color = black] (0,-8) rectangle  (6,-9);

 \node[] at (0,-6.8) 
    {$0$};
 \node[] at (1,-6.8) 
    {$1$};
 \node[] at (2,-6.8) 
    {$2$};
 \node[] at (3,-6.8) 
    {$3$};
 \node[] at (4,-6.8) 
    {$4$};
 \node[] at (5,-6.8) 
    {$5$};
 \node[] at (6,-6.8) 
    {$6$};
 \node[] at (0,-6.4) 
    {time $\rightarrow$};

 \node[] at (-1,-7.5) 
    {job $1$};
 \node[] at (-1,-8.5) 
    {job $2$};

 \node[text width=6cm] at (3,-9.5) 
    {Third proposed schedule};

 \node[text width=3cm] at (8,-7.5) 
    {$C_{\max} = \frac{51}{10}$,};

 \node[text width=3cm] at (8,-8.5) 
    {$\sum_{j \in J}C_j = \frac{101}{10}$};


\fill[color = gray] (0,-12) rectangle  (5,-11.6);
	

\fill[color = gray] (0,-12.6) rectangle  (5,-13);

\draw[color = red] (0,-13) -- (5,-12);
\draw[color = red] (0,-12) -- (5,-11);

\draw[color = black] (0,-12) rectangle  (6,-11);
\draw[color = black] (0,-12) rectangle  (6,-13);

 \node[] at (0,-10.8) 
    {$0$};
 \node[] at (1,-10.8) 
    {$1$};
 \node[] at (2,-10.8) 
    {$2$};
 \node[] at (3,-10.8) 
    {$3$};
 \node[] at (4,-10.8) 
    {$4$};
 \node[] at (5,-10.8) 
    {$5$};
 \node[] at (6,-10.8) 
    {$6$};
 \node[] at (0,-10.4) 
    {time $\rightarrow$};

 \node[] at (-1,-11.5) 
    {job $1$};
 \node[] at (-1,-12.5) 
    {job $2$};

 \node[text width=6cm] at (3,-13.5) 
    {Schedule with fractional loads};

 \node[text width=3cm] at (8,-11.5) 
    {$C_{\max} = 5$\,};

 \node[text width=3cm] at (8,-12.5) 
    {$\sum_{j \in J}C_j = 10$};

\end{tikzpicture}
    \caption{An illustration of how preempting jobs after smaller increments in time can lead to schedules with smaller completion times. The gray areas represent periods during which a job is processed, and the red line indicates the temperature of a job.}
    \label{fig:example}
\end{figure}

Let us first focus on a naive schedule one could construct for the two jobs. One can take the aforementioned schedule for job 1, and then process job 2 at times when the machine is not processing job 1. This leads to a completion time of~$6$ for job 2, see the first proposed schedule in Figure~\ref{fig:example}. This schedule has~$C_{\max}=6$ and~$\sum_{j \in J}C_j=11$. 

However, one can actually decide to stop processing a job prior to the job hitting the maximum temperature~$T=1$. For instance, the second proposed schedule in Figure~\ref{fig:example} considers time intervals of length~$\frac{1}{5}$ and processes job 1 alternately in these time intervals. Job 2 is processed in a similar way, but with a delay of~$\frac{1}{5}$. The resulting schedule has~$C_{\max} = \frac{26}{5}$ and~$\sum_{j \in J}C_j = \frac{102}{10}$. And when we make the time intervals even smaller, we can get a schedule with~$C_{\max} = \frac{51}{10}$ and~$\sum_{j \in J}C_j = \frac{101}{10}$, see the third proposed schedule in Figure~\ref{fig:example}. 

The final schedule in Figure~\ref{fig:example} illustrates that by using increasingly small time intervals we may  obtain a schedule with~$C_{\max}\leq 5+\epsilon$ and~$\sum_{j \in J}C_j\leq 10+\epsilon$  for any~$\epsilon>0$. This schedule has fractional workloads representing the proportion of time each job is fully loaded on the machine in a sufficiently fine alternating schedule. However, no schedule in which the total number of job pre-emptions is finite can attain the value of this fractional schedule, because in any such schedule there is a job that is not processed at all until time~$\epsilon>0$. 
This example shows that to approach the infima of interest, we may need to take into account infinitesimal amounts of processing time, and use an arbitrary number of preemptions in our schedules.  At the same time the example suggests that we may characterise the value of the infimum, as well as a scheme for creating schedules that approach the infimum, by means of a schedule with fractional workloads.
In the next section, we show that minimizing the sum of completion times / makespan over a certain class of fractional schedules indeed is a viable approach to our problems.

\section{Preparatory Insights}
\label{sec:prepinsight}
In this section, we broaden the class of schedules considered so far to include schedules with fractional workloads. The schedules of our original problem definition will be named {\em natural} schedules, and we introduce {\em simple} and {\em normal} schedules. We show that the infimum of the sum of completion times (or the makespan) over either of these three classes is the same, and that the infimum over normal schedules is attained.

Formally, a {\em schedule}~$S=(S_j: j\in J)$ for a set of jobs~$J$ consists of a {\em load function}
$$S_j:\R_+\rightarrow \R_+$$ for each job~$j\in J$, indicating the fractional load of job~$j$ at time~$t$. 
To define the feasibility conditions and the objective of our scheduling problem, we must assume that these load functions are integrable. Note that this will be the case for the natural, simple, and normal schedules formally defined below.

A schedule is {\em manageable} if the total load never exceeds~$1$, that is, if~$$\sum_{j\in J} S_j(t)\leq 1$$ for all~$t\geq 0$.  
A schedule~$S$ will uniquely determine, for each job~$j\in J$, a {\em temperature}~$T_j(t)$ so that 
\begin{equation}
\label{eq:derivtemp} 
\begin{aligned}\frac{\partial}{\partial t}&T_j(t)=\\ & \hspace*{-0.83cm}\left\{\begin{array}{ll}
\alpha(1-S_j(t))+\beta S_j(t)&\text{ if }T_j(t)> 0\\
\max\{0,\alpha(1-S_j(t))+\beta S_j(t)\}&\text{ if }T_j(t)= 0.
\end{array}\right.
\end{aligned}\end{equation}
We say that the schedule~$S$ is {\em feasible} if it is manageable and~$T_j(t)\leq 1$ for all~$j\in J$ and~$t\geq 0$. 
The total work on job~$j$ at time~$t$ is captured by
$$W_j(t):=\int_0^t S_j(t) dt.$$
Given a schedule~$S$ and processing times~$p_j$ for each job~$j\in J$, the 
{\em completion time } of job~$j\in J$ is determined as 
$$C_j:=\min\{t\in \R_+: p_j\leq W_j(t)\},$$
and we write~$C(S)=(C_j: j\in J)$.

We say that a schedule is {\em natural} if~$S_j(t)\in \{0,1\}$ for all~$j\in J, t\geq 0$, and~$\{t\in \R_+: f(t)=1\}$ is a union of half-open intervals~$[t, t')$. 
That is, in a natural schedule a job is fully loaded on the machine or not at all. In its original form, our main problem is to minimize the makespan, resp. the sum of completion times, over all natural schedules. 

We will call a schedule~$S$ {\em simple} if there are~$0=t_0<t_1<\cdots<t_m$ so that 
$$S_j(t)=S_j(t')\text{ for all }t, t'\in [t_{i-1},t_i)$$
 for all~$j\in J$ and~$i=1,\ldots, m$. We then say that the simplicity of~$S$ is {\em witnessed } by~$\{t_1,\ldots, t_m\}$. A {\em normal} schedule is a simple schedule whose simplicity is witnessed by~$\{C_j(S): j\in J\}$. That is, in a normal schedule~$S$ the fractional loads are constant between the successive completion of two jobs.

\ignore{
Natural schedules are unwieldy objects. Even if we insist that the load functions be integrable and that jobs are loaded on the machine for intervals of non-zero length, a natural schedule may involve an indefinite number of times where the machine preempts, i.e. stops the processsing of the current jobs and loads another job. As a result, there is no upper bound on the length of a description of a natural schedule in terms of the number of jobs, say. Testing feasibility may already involve the evaluation of infinite sums, let alone optimizing over all natural schedules. When minimizing the sum of completion times over all natural schedules, the infimum may not be attained. E.g. in the example with 2 jobs illustrated in Figure \ref{fig:example}, each job takes at least 5 units of time to complete from the moment the job is first loaded on the machine. Since the machine can only load one of the jobs at time~$t=0$, there will necessarily be a slight lag before the other job can start, so that that the sum of completion times will always be slightly more than~$5+5$ in a natural schedule. 
}

Natural schedules are unwieldy objects. In a natural schedule, each load function~$S_j$ is determined by the set of times~$t$ where~$S_j$ toggles on or off, but there is no a priori bound on the cardinality of this set of toggle points. This limits the efficiency of any algorithm that handles explicit descriptions of such schedules.

By contrast, a normal schedule for~$n$ jobs is fully determined by~$O(n^2)$ parameters. Assuming for simplicity that~$J=\{1,\ldots, n\}$, let~$C\in \R_+^n$ be such that~$C_1\leq C_2\leq \cdots\leq C_n$ and let~$W\in \R_+^{n\times n}$ be such that
\begin{equation}\label{eq:nondec}W_{1,j}\leq\cdots\leq W_{n,j}\text { for } j\in J.\end{equation}
Then the unique normal schedule~$S$ so that~$C(S)_{j}=C_j$ and~$W_j(C_i)=W_{i,j}$ is given by 
$$S_j(t)=\frac{W_{i,j}-W_{i-1,j}}{C_{i}-C_{i-1}}\text{ for }t\in [C_{i-1}, C_i)$$
for~$i \in J$, taking~$C_0:=0$. 

The feasibility of the normal schedule~$S$ can be expressed in terms of linear inequalities on~$C, W$ as follows. 
Ensuring that job~$j$ completes at time~$C_j$ in schedule $S$ amounts to 
\begin{equation}\label{eq:complete}W_{i,j}=p_j\text{ for }i,j\in J, i\geq j.\end{equation}
The schedule~$S$ is manageable if and only if 
\begin{equation}\label{eq:man} \begin{split}&\sum_{j\in J} W_{1,j}\leq C_1\\&\text{and for }i\in J \setminus 1:\\
&\sum_{j\in J} (W_{i,j}-W_{i-1,j})\leq C_{i}-C_{i-1}.\end{split}\end{equation}
Since the fractional load of job~$j$ is constant throughout the interval~$[C_{i-1}, C_i)$, we have~$T_j(t)\leq \max\{T_j(C_{i-1}), T_j(C_i)\}$ for all~$t\in [C_{i-1}, C_i)$. So the feasibility condition that~$T_j(t)\leq 1$ for all~$t$ reduces to~$T_j(C_i)\leq 1$ for~$i=1,\ldots, n$. The temperature of job~$j$ at completion time of job~$i$ is determined by~$T_j(C_i)=\max\{0, T_j(C_{i-1})+\alpha(C_i-C_{i-1}-W_{i,j}-W_{i-1,j})+\beta (W_{i,j}-W_{i-1,j})$. Therefore, the schedule~$S$ has a feasible temperature profile if and only if there exists a~$T\in \R_+^{n\times n}$ so that 
\begin{equation}\label{eq:temp}\alpha C_1+(\beta-\alpha) W_{1,j}\leq T_{1,j}\text{ for } ~j\in J,\end{equation}
\begin{equation}\label{eq:temp2}
\begin{aligned}\alpha(C_i-C_{i-1})+(\beta-\alpha) (W_{i,j}-W_{i-1,j})\\\leq T_{i,j}-T_{i-1,j}\text{ for }i, j\in J, i\neq 1,\end{aligned}\end{equation}
and 
\begin{equation}\label{eq:temp3}T_{i,j}\leq 1\text{ for all } ~i, j\in J.\end{equation}
Summarizing the above discussion, we have:
\begin{lemma} \label{lem:norm} Let~$J=\{1,\ldots n\}$ and let~$C\in \R^n_+, W\in \R_+^{n\times n}$.
Suppose that 
$$C_1\leq C_2\cdots\leq C_n.$$
Then the following are equivalent.
\begin{enumerate}
\item There is a normal schedule~$S$ so that~$$C(S)_j=C_j\text{ and }W_j(C_i)=W_{i,j}\text{ for all }i,j\in J.$$
\item There exists a~$T\in \R^{n\times n}_+$ so that~$C,W,T$ satisfies \eqref{eq:nondec},
\eqref{eq:complete}, \eqref{eq:man}, \eqref{eq:temp}, \eqref{eq:temp2}, and \eqref{eq:temp3}.
\end{enumerate}
\end{lemma}

Lemma~\ref{lem:norm} allows us to prove the main result of this section, namely that optimizing over all natural schedules is equivalent to optimizing over all normal schedules, given that the objective is a continuous function of the completion times.
\begin{theorem}\label{thm:natural=normal}Let~$f:\R^n\rightarrow \R$ be continuous. Then 
\small{
\begin{align}\label{eq:nat}&\inf\left\{f(C(S)): S\text{ a feasible natural schedule}\right\}=\\&
\label{eq:simple}\inf\left\{f(C(S)): S\text{ a feasible simple schedule}\right\}= \\&
\label{eq:nor}\inf\left\{f(C(S)): S\text{ a feasible normal schedule}\right\}. 
\end{align}
}
\end{theorem}
\proof 
$\eqref{eq:nat}\geq \eqref{eq:nor}$: Let~$S$ be a feasible natural schedule. Without loss of generality, we may assume that~$J=\{1,\dots, n\}$ and~$C_1\leq\cdots\leq C_n$ for~$C=C(S)\in \R^n_+$. Let   
$W, T\in \R^{n\times n}_+$ be given by~$$W_{i,j}:=W_j(C_i)\mbox{ and } T_{i,j}:=T_j(C_i).$$ It is straightforward that the triple~$C, W, T$ satisfies \eqref{eq:nondec}, \eqref{eq:complete}, \eqref{eq:man}, \eqref{eq:temp}, \eqref{eq:temp2}, and  \eqref{eq:temp3}.
By Lemma \ref{lem:norm}, there is a normal schedule~$S'$ so that~$C(S')=C$. Then~$\eqref{eq:nor}\leq f(C(S'))\leq f(C(S))$, as required. 

$\eqref{eq:nor}\geq \eqref{eq:simple}$: A normal schedule is also a simple schedule, hence any feasible solution for \eqref{eq:nor} is also feasible for \eqref{eq:simple}.

$\eqref{eq:simple}\geq \eqref{eq:nat}$: Let~$S$ be a feasible simple schedule. Say, there are~$0=t_0<t_1<\cdots<t_m$ so that for~$i=1,\ldots, m$ we have~$S_j(t)=S_j(t')$ for~$t,t'\in [t_{i-1}, t_i)$.

Fix a~$\gamma>1$, and let~$S^{\gamma}$ be the schedule  specified by
$$S^{\gamma}_j(\gamma t)=\frac{S_j(t)}{\gamma}\qquad\text{ for all }t\in \R_+.$$
In the schedule~$S^\gamma$, the total work on job~$j$ during an interval~$[\gamma t_{i-1}, \gamma t_i)$ equals the work on~$j$ during~$[t_{i-1}, t_i)$ in the original schedule~$S$, and therefore~$C(S^\gamma)=\gamma C( S)$. The total workload of the processor reduces by a factor~$\gamma$ throughout, and hence the schedule~$S^\gamma$ will be manageable. In~$S^\gamma$, we take more time to accomplish the same amount of work on~$j$ in corresponding intervals, hence there is a~$\kappa$ so that~$T^{\gamma}(t)\leq\kappa<1$ throughout for the temperature function of~$S^\gamma$.

Let~$S^{\gamma,k}$ be the natural schedule that arises from~$S^{\gamma}$ by discretizing and time-slicing each interval~$[\gamma t_{i-1}, \gamma t_i)$. 
That is, in~$S^{\gamma,k}$ the interval~$[\gamma t_{i-1}, \gamma t_i)$ is subdivided in~$k$ intervals of equal length which repeat the same natural micro-schedule. For each job~$j$, that micro-schedule has an interval of length~$w_j\gamma(x_i-x_{i-1})/k$ where~$j$ is fully active and all other jobs are inactive, where~$w_j:=S^{\gamma}(t)$ for any~$t\in [\gamma t_{i-1}, \gamma t_i)$. Thus the length of the interval for~$j$ in the micro-schedule is proportional to the constant fractional workload~$w_j$ of~$j$ during~$[\gamma t_{i-1}, \gamma t_i)$.
The micro-schedule will necessarily have times where no job is active, since~$\sum_{j \in J} w_j\leq 1/\gamma<1$. 
It is evident that~$S^{\gamma,k}$ is manageable, and that~$S^{\gamma,k}$ and~$S^\gamma$ do the same total work on each job~$j$ in each interval~$[\gamma x_{i-1}, \gamma x_i)$.

For any~$\epsilon>0$, there is a sufficiently high~$k$ so that~$$|T^{\gamma,k}(t)-T^{\gamma}(t)|<\epsilon$$
for all~$t\in \R_+$. Since~$T^{\gamma}(t)\leq \kappa<1$, it follows that there is a large enough~$k(\gamma)$ so that~$T^{\gamma,k(\gamma)}(t)\leq 1$ for all~$t\in \R_+$. In other words,~$S^{\gamma,k(\gamma)}$ is feasible. Finally, we may assume that~$k(\gamma)\rightarrow\infty$ as~$\gamma\downarrow 1$, so that~$
C(S^{\gamma,k(\gamma)})\rightarrow C(S^\gamma)=\gamma C(S)$ as~$\gamma\downarrow 1$.
It follows that for any feasible simple schedule~$S$, we have {\small{
\[\eqref{eq:nat}\leq \lim_{\gamma\downarrow 1} f(C(S^{\gamma,k(\gamma)}))=\lim_{\gamma\downarrow 1} f(\gamma C(S))= f(C(S))\]}}
as~$f$ is continuous. Then~$\eqref{eq:nat}\leq \eqref{eq:simple}$, as required.
\endproof

Applying the theorem to~$f(S)=\sum_j C_j(S)$ (resp.~$f(S)=\max_j C(S)_j$), we find that minimizing the sum of completion times (resp. the makespan) over all natural schedules equals minimizing over all normal schedules. At the same time, Lemma~\ref{lem:norm} implies that minimizing~$f(S)$
over all normal schedules~$S$ so that~$C_1(S)\leq \cdots\leq C_n(S)$ amounts to solving the following linear program (LP).
\begin{strip}
\begin{equation}
\label{LP}
\begin{array}{lll}\min & \sum_{j\in J} C_j\\
\\
&W_{1,j}\leq\cdots\leq W_{n,j}\\
&W_{i,j}=p_j &~i,j\in J, i\geq j\\
\\
&\sum_{j\in J} W_{1,j}\leq C_1\\
&\sum_{j\in J} (W_{i,j}-W_{i-1,j})\leq C_{i}-C_{i-1}&~i\in J, i\neq 1\\
\\\
&\alpha C_1+(\beta-\alpha) W_{1,j}\leq T_{1j}& ~j\in J\\
&\alpha(C_i-C_{i-1})+(\beta-\alpha) (W_{ij}-W_{i-1,j})\leq T_{i,j}-T_{i-1,j} &~i, j\in J, i\neq 1\\
&T_{i,j}\leq 1& ~i, j\in J\\
\\
&C\in \R_+^n, ~W, T\in \R^{n\times n}_+
\end{array}
\end{equation} 
\end{strip}
Since this LP has ~$O(n^2)$ variables and~$O(n^2)$ linear constraints with coefficients whose bit-lengths are bounded by the bit-length of the input parameters~$(\alpha, \beta, p)$, solving  \eqref{LP} takes polynomial time.

For each linear ordering~$<$ of the jobs~$J$, we can evidently draw up a similar LP to compute a best schedule that finishes the jobs in the order specified by~$<$. Indeed: 
\begin{equation}\label{eq:order}
\begin{aligned} &
z^<:=\\& \hspace*{-0.72cm}\min\left\{\sum_{j\in J} C_j(S): ~\begin{aligned}S\text{ is feasible and normal, }\\ i<j \Rightarrow C_i(S)\leq C_j(S)\end{aligned}\right\} \end{aligned}\end{equation} 
This immediately implies that  the decision problem associated with minimizing the sum of completion times over all normal schedules is in NP. 
Indeed, if~$\sum_{j\in J} C_j(S)\leq c$ for some feasible normal schedule~$S$, then any linear ordering~$<$ of~$J$ so that~$i<j \Rightarrow C_i(S)\leq C_j(S)$ certifies this fact in polynomial time, by checking that ~$z^<\leq c$.

Computing the overall optimal normal schedule
\begin{equation}\label{eq:normal} \hspace*{-0.4cm}\min\left\{\sum_{j\in J} C_j(S): ~S\text{ is feasible 
 and normal}\right\}\end{equation} 
in polynomial time is still another matter. Clearly \eqref{eq:normal} equals 
\begin{equation}\label{eq:orders}\min\{z^<: <\text{ a linear order of }J\}.\end{equation}
But even if evaluating any~$z^<$ takes polynomial time, this will not suffice if we have to evaluate~$z^<$ for all~$n!$ linear orders~$<$ of~$J$. 
To overcome this, we will identify an order of completion~$<$ that attains the minimum in \eqref{eq:orders}. We establish this result in Section~\ref{sec:sumofcompletiontimes}.

\subsection{The example continued}
We revisit the example of Section \ref{ssec:example}, with~$J=\{1,2\}$,~$p_1 = p_2 = 2, \alpha = -\frac13$,~$\beta = 1$. Due to the symmetry between the two jobs, we may assume that~$C_1\leq C_2$ in some normal schedule minimizing the sum of completion times. Formulation~(\ref{LP}) applied to the example gives the following LP.
\begin{equation}
\nonumber
\begin{array}{lll}\min & C_1+C_2\\
\\
&W_{1,1}= W_{2,1}=2\\
&W_{1,2}\leq W_{2,2}=2 \\
\\
&W_{1,1}+W_{1,2}\leq C_1\\
&W_{2,1}-W_{1,1}+W_{2,2}-W_{1,2}\leq C_2-C_1\\
\\
&-\frac13 C_1+\frac43 W_{1,1}\leq T_{1,1}\\
&-\frac13 C_1+\frac43 W_{1,2}\leq T_{1,2}\\
&-\frac13(C_2-C_1)+\frac43 (W_{2,1}-W_{1,1})\leq T_{2,1}-T_{1,1}\\
&-\frac13(C_2-C_1)+\frac43 (W_{2,2}-W_{1,2})\leq T_{2,2}-T_{1,2}\\
&T_{1,1}\leq 1\\
&T_{1,2}\leq 1\\ &T_{2,1}\leq 1\\ &T_{2,2}\leq 1& \\
\\
&C\in \R_+^2, ~W, T\in \R^{2\times 2}_+
\end{array}
\end{equation}
An optimal solution to this LP is~$C_1=C_2=5$,~$W_{i,j}=2$,~$T_{i,j}=1$ for all~$i,j$, leading to load functions~$S_1(t) = S_2(t) = \frac25$ for~$0 \leq t \leq 5$; see Figure \ref{fig:example}.

\section{Minimizing the sum of completion times}
\label{sec:sumofcompletiontimes}

In this section, we show how to find a schedule minimizing \(\csum\) in polynomial time.
A first and crucial step to find an optimal schedule for \(n\) jobs with respect to \(\csum\), is to establish the order in which the jobs finish. 

\begin{theorem}
\label{thm:fixedvolgorde}
Given an instance of~$Q1^{sum}$, there is an optimal normal schedule where the jobs are completed in order from smallest to largest processing times. 
\end{theorem}
 With this proven, the LP \eqref{LP} from Section~\ref{sec:prepinsight} solves the problem. Clearly, as this LP can be solved in polynomial time, we can find an optimal schedule in polynomial time. 
 \begin{corollary}
~$Q1^{sum}$ is solvable in polynomial time.
\end{corollary}
The bulk of the proof of Theorem \ref{thm:fixedvolgorde} will be the verification of the following key lemma.
\begin{lemma}\label{lem:key} Let~$S$ be a feasible normal schedule, and  
suppose that~$i,j\in J$ are such that~$p_i<p_j$ and~$C_i(S)>C_j(S)$. Then there exists a feasible normal schedule~$S'$ so that~$$C_k(S')=\left\{\begin{array}{ll}C_k(S)&\text{if }k\not\in\{i,j\},\\
C_i(S)&\text{if }k=j,\\
C_j(S)&\text{if }k=i.\\
\end{array}\right.$$
\end{lemma}
Having this lemma, the proof of the theorem is straightforward.
\proof[Proof of Theorem~\ref{thm:fixedvolgorde}]
Let~$J=\{1,\ldots, n\}$. We may assume that~$p_1\leq \cdots\leq p_n$.

Choose an optimal schedule~$S$ so that the number of pairs~$(i,j)$ that complete in the wrong order,
$$m(S):=\#\{(i,j)\in J\times J: p_i<p_j, C_i(S)>C_j(S)\},$$
is as small as possible. If there exists an~$i$
so that~$C_i(S)>C_{i+1}(S)$, then by applying Lemma \ref{lem:key} to~$S$ and~$(i,j)=(i,i+1)$, we obtain a schedule~$S'$  so that~$m(S')=m(S)-1<m(S)$. This would contradict the choice of~$S$. It follows that ~$C_1(S)\leq \cdots\leq C_n(S)$, as required.
\endproof

\proof[Proof of the Lemma~\ref{lem:key}]
For each job~$k\in J$, let~$T_k$ resp.~$W_k$ denote the temperature resp. work functions for job~$k$ associated with the schedule~$S$. 

We will first establish the lemma under the following assumption on~$S$:
 \begin{equation}\label{eq:tempnot1} T_k(t)<1\text{ for all } t\in\R_+\end{equation}
The general case of the lemma will follow a continuity argument based on this special case.

To prove the lemma, we will create the schedule \(S'\) by finding a specific point in time \(\tequal\) where both jobs \(i,j\) are equally far away from completion. We show that we can swap the fractional load of jobs~$i$ and~$j$ from \(\tequal\) on wards. This creates a schedule~$S'$ that switches the completion times of~$i$ and~$j$, as both jobs are equally far away from completion at~$\tequal$. The new schedule~$S'$ is evidently manageable as the total load of the machine is unaffected. To make sure that the new schedule is feasible, it remains to show that for the temperatures induced by~$S'$ we have~$T'_k(t)\leq 1$. As the temperature functions~$T'_k$ are unchanged for~$k\not \in\{i,j\}$, we need to do so specifically for~$T'_i(t)$ and~$T_j'(t)$ .

We note that the schedule~$S'$ constructed by switching this way is simple but may not be a normal schedule, since due to the switching at~$\tequal$ the simplicity of~$S'$ is witnessed by~$\{\tequal\}\cup\{C_k(S): k\in J\}$. One may obtain a feasible normal schedule from~$S'$ by averaging each load function~$S'_i$ and~$S'_j$ over the interval between the two completion times surrounding~$\tequal$.

For any job~$k \in J$, we let \(\lefttodo_k(t)\) indicate how much of job \(k\) still needs to be processed at time \(t\), i.e.~$$\lefttodo_k(t) := p_k-W_k(t) = p_k - \int_{0}^{t} S_k(t') dt'.$$ 
We also define  \[\T = \{t : \lefttodo_i(t) = \lefttodo_j(t)\}\] as the times where jobs~$i$ and~$j$ are equally far away from completion. 
We will heavily use the set~$\mathcal{T}$, and will use the following observation a multitude of times.

\begin{observation} \label{obs:1}
    If~$\lefttodo_i(t) \ge \lefttodo_j(t)$ and~$\lefttodo_i(t') \le \lefttodo_j(t')$ for any~$t$ and~$t'$, then there exists~$\bar{t}$ in between~$t$ and~$t'$ with~$\bar{t} \in \mathcal{T}$.  
\end{observation}

Indeed, since the functions~$\lefttodo_i(t)$ and~$\lefttodo_j(t)$ are continuous functions, at one point they then should be equal. 
We distinguish three types of points in \(\T\):
\begin{itemize}
\item \(\T^> = \{t \in \T \colon T_i(t) > T_j(t)\}\),
\item \(\T^= = \{t \in \T \colon T_i(t) = T_j(t)\}\),
\item \(\T^< = \{t \in \T \colon T_i(t) < T_j(t)\}\).
\end{itemize}
Using Observation~\ref{obs:1} we see~$\T \not = \emptyset$ as~$p_i = \lefttodo_i(0) < \lefttodo_j(0) = p_j$ and because~$C_i > C_j$ we have~$\lefttodo_i(C_j) > \lefttodo_j(C_j) = 0$.

Note that if there exists~$t \in \T^=$, then we can switch the processing of jobs~$i$ and~$j$ from~$t$ on wards, and since the temperatures of both jobs were equal at~$t$, this directly proofs that this is a feasible schedule. Therefore, we will assume from now on that~$\T^= = \emptyset$. 

We will use the following claim often to argue about our schedules:

\begin{claim}\label{claim:workdone}
For job~$k\in J$ and times~$t_1 \le t_2$ we have
    \[\lefttodo_k(t_1) - \lefttodo_k(t_2) \le \frac{T_k(t_2) - T_k(t_1) - \alpha (t_2-t_1)}{\beta-\alpha}.\]
    Moreover, if~$T_k(t) >0$ for all~$t \in (t_1,t_2)$, then
    \[\lefttodo_k(t_1) - \lefttodo_k(t_2) = \frac{T_k(t_2) - T_k(t_1) - \alpha (t_2-t_1)}{\beta-\alpha}.\]
\end{claim}
\begin{claimproof}
    According to (\ref{eq:derivtemp}) we have that~$\frac{\partial}{\partial t} T_k(t) \ge \alpha(1- S_k(t)) + \beta S_k(t)$. Hence, we get \[S_k(t) \le \frac{\frac{\partial}{\partial t}T_k(t) - \alpha}{\beta - \alpha}.\] By taking the integral on both sides from~$t_1$ to~$t_2$, together with noting that~$\lefttodo_k(t_1) - \lefttodo_k(t_2) = \int_{t_1}^{t_2}S_k(t) \,dt$, we get the requested inequality. Moreover, note that 
   ~$\frac{\partial}{\partial t} T_k(t) = \alpha(1- S_k(t)) + \beta S_k(t)$ if~$T_k(t) >0$ (see (\ref{eq:derivtemp})). Hence, if~$T_k(t) >0$ for all~$t \in (t_1,t_2)$ we get that all inequalities becomes equalities.
\end{claimproof}

We want to define~$\tequal \in \T^<$ as the moment to switch the processing of jobs~$i$ and~$j$. For this reason we first prove that this set is nonempty.
\begin{claim}\label{claim:2} $\T^< \not = \emptyset$.
\end{claim}
\begin{claimproof}
We prove this by contradiction. Take~$t^\star = \min \T$. As we assumed~$\T^= = \emptyset$ and~$\T^< = \emptyset$, we find that~$t^\star \in \T^>$. Let  
\[t^= = \max \{t: 0 \le t < t^\star: T_i(t) = T_j(t)\}.\]
Note that this~$t^=$ is well-defined as~$T_i(0) = T_j(0) = 0$. From the definition of~$t^=$ we find that~$T_i(t) > T_j(t)$ for all~$t \in (t^=, t^\star]$. In particular, we find~$T_i(t) > 0$ for all~$t \in (t^=, t^\star]$.
Therefore we get using Claim~\ref{claim:workdone}:
\begin{align*}
    \lefttodo_j(t^=) - \lefttodo_j(t^\star) &\le  \frac{T_j(t^\star) - T_j(t^=) - \alpha(t^\star - t^=)}{\beta - \alpha} \\
& = \frac{T_j(t^\star) - T_i(t^=) - \alpha(t^\star - t^=)}{\beta - \alpha}\\
&\le \frac{T_i(t^\star) - T_i(t^=) - \alpha(t^\star - t^=)}{\beta - \alpha} \\
&= \lefttodo_i(t^=) - \lefttodo_i(t^\star).
\end{align*}

In other words, more work has been done on~$i$ than on~$j$ during the interval~$[t^=,t^\star]$. Combining this fact with 
$\lefttodo_i(t^\star) =\lefttodo_j(t^\star)$ as~$t^\star \in \T$, we get that~$\lefttodo_i(t^=) > \lefttodo_j(t^=)$. 
However, since~$\lefttodo_i(0) = p_i < p_j =\lefttodo_j(0)$, we find (using Observation~\ref{obs:1}) that there exists~$\bar{t} < t^\star$ such that~$\bar{t} \in \T$, a contradiction to the choice of~$t^\star = \min \T$.
\end{claimproof}

We define
\begin{align*}
\tequal =  \max \T^< 
\end{align*}
as the moment at which will we switch the processing of jobs~$i$ and~$j$.
Using \(\tequal\), we construct a schedule \(S'\) where~$S'_k(t) = S_k(t)$ for all~$t$ and all jobs~$k \in J$, with the exception that~$S'_i(t) = S_j(t)$ and~$S'_j(t) = S_i(t)$ for all~$t \ge \tequal$. We use~$T'_k(t)$ to refer to the temperature of job~$k \in J~$ at time~$t$ in schedule~$S'$.

\paragraph*{Proving that schedule~$S'$ is feasible.} 

Since \(\tequal \in \T^<\), we immediately see that job~$i$ does not violate the temperature constraint in~$S'$ as \(T_i(\tequal) < T_j(\tequal)\) and~$i$ is processed in~$S'$ in the same way as~$j$ is processed in~$S$. As all other jobs in~$S'$ are the same way processed as in~$S$, we only need to show that job~$j$ does not violate any temperature constraints in~$S'$, that is~$T_j'(t) \leq 1$ for all~$t > \tequal$. 

Assume for the purpose of contradiction that~$j$ is too hot directly after time \(t' > \tequal\), i.e.~$T'_j(t') = 1$ where~$t'$ is the first moment where this happens. Note that we can't have~$t' = \tequal$ because~$T'_j(\tequal) = T_j(\tequal) < 1$ by \eqref{eq:tempnot1}. 

\begin{claim} \label{claim:t<exists} \label{claim:3}
    There exists \(t^> = \min\{t \in \T^> : t > \tequal\}\).
\end{claim}
\begin{claimproof}
At time \(t'\) we have \(T'_j(t') = 1\). If we would have~$T'_j(\bar t)=0$ for any~$\bar t \in [\tequal,t']$, we would find that~$T'_j(\bar t) \le T_i(\bar t)$ and since~$j$ is processed in~$S'$ in the same way as~$i$ is processed in~$S$ we find that~$T'_j(t) \le T_i(t) \le 1$ for all~$t \ge \bar t$. In particular, it would contradict that job~$j$ is too hot directly after time~$t' \ge \bar t$. Hence, we can assume that~$T'_j(t) > 0$ for all~$t \in [\tequal,t']$. 

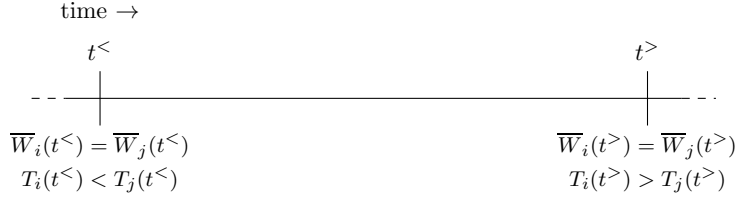
\begin{figure*}[h!]
    \centering
    \begin{tikzpicture} [scale = .9, transform shape,
]

\draw[dashed] (-1,0) -- (-.5,0);
\draw[color = black] (-0.5,0) -- (8.5,0);
\draw[dashed] (8.5,0) -- (9,0);

 \node[] at (0,1.3) 
    {time $\rightarrow$};
 \node[] at (0,.7) {$\tequal$};
 
 \node[] at (0,-.7) {\small $\lefttodo_i(\tequal) = \lefttodo_j(\tequal) $};
\node[] at (0,-1.2) {\small $T_i(\tequal) < T_j(\tequal) $};
\draw[] (0,-.4) -- (0,.4);

\node[] at (8,.7) {$t^>$};
\draw[] (8,-.4) -- (8,.4);
\node[] at (8,-.7) {\small $\lefttodo_i(t^>) = \lefttodo_j(t^>) $};
\node[] at (8,-1.2) {\small $T_i(t^>) > T_j(t^>) $};


\end{tikzpicture}
    \caption{Claim~\ref{claim:2} ensures the existence of $t^< \in \T^<$. The properties of $t^<$ follow directly from $t^< \in \T^<$. Similarly, Claim~\ref{claim:3} ensures the existence of $t^> = \min \{t \in \T^>: t>t^<\}$. As $t^> \in \T^>$ the listed properties follow.
    }
    \label{fig:timeline1}
\end{figure*}

As a consequence, we can exactly compute how much work on \(j\) has been done in~$S'$ in the period~$(\tequal, t' )$ using Claim~\ref{claim:workdone}, namely:
\begin{align*}
&\frac{T'_j(t') - T'_j(\tequal) - \alpha(t'- \tequal)}{\beta - \alpha} \\&\qquad>  \frac{T_j(t') - T_j(\tequal) - \alpha(t'- \tequal)}{\beta - \alpha} \\
&\qquad\ge \lefttodo_j(t^<) - \lefttodo_j(t'),
\end{align*}
where in the first inequality we use~$T'_j(t') = 1 > T_j(t')$ and~$T_j(\tequal) = T'_j(\tequal)$. Now note that~$\int_{t^<}^{t'}S'_j(t)\,dt = \int_{t^<}^{t'}S_i(t)\,dt$ (i.e. the work done between~$t^<$ and~$t'$ is the same for~$i$ in~$S$ as it is for~$j$ in~$S'$,) and so we can conclude that 
\[\lefttodo_i(t^<) - \lefttodo_i(t') > \lefttodo_j(t^<) - \lefttodo_j(t').\]
Since~$\lefttodo_i(t^<) =\lefttodo_j(t^<) ~$ as~$t^< \in \T$, we get~$\lefttodo_i(t') < \lefttodo_j(t')$.
However, since we also have~$\lefttodo_j(C_j) = 0 < \lefttodo_i(C_j)$ as job~$i$ is finished after job~$j$ in~$S$, we find using Observation~\ref{obs:1} that there must exist a~$t > t'~$ such that~$\lefttodo_i(t) = \lefttodo_j(t)$. This implies that~$\{t \in \T : t > \tequal\}$ is nonempty and since~$\tequal$ as chosen as maximal~$t \in \T^<$, we find that any~$t \in \T$ after~$\tequal$ must be part of~$\T^>$. 

So, we proved that ~$t^> = \min \{t \in \T^>:  t> \tequal\}$
is well-defined.
\end{claimproof}

See Figure~\ref{fig:timeline1} for an illustration. Note that we have the following observation as a consequence of the definitions of~$\tequal$ and~$t^>$.

\begin{observation}\label{obs2}
For all~$t \in (\tequal,t^>)$ we have~$t \not \in \T$.
\end{observation}
This observation will be used to later.

\begin{claim}\label{claim:4}
There exists \[t^= =\max \{t \in (\tequal,t^>): T_i(t) = T_j(t)\}.\]    
\end{claim}
\begin{claimproof}
    Note that~$T_i(\tequal) < T_j(\tequal)$ and~$T_i(t^>) > T_j(t^>)$ as~$\tequal \in \T^<$ and~$t^> \in \T^>$. As the temperatures are continuous functions, there must be time where they are equal. Hence
$t^= = \max \{t \in (\tequal,t^>): T_i(t) = T_j(t)\}$
is well-defined.
\end{claimproof}

\begin{claim}
     \label{claim:5}
     We have~$\lefttodo_i(t^=) > \lefttodo_j(t^=)$.
\end{claim}

\begin{claimproof}
    For this, we actually prove that there is a~$t_0 \in (\tequal,t^>)$ with~$\lefttodo_i(t_0) > \lefttodo_j(t_0)$. This is indeed sufficient:~$\lefttodo_i(t^=) \le \lefttodo_j(t^=)$ combined with this would imply (using Observation~\ref{obs:1}) that there exists a~$t \in (t^<, t^>)$ between~$t_0$ and~$t^=$ such that~$t \in \T$. However, this would contradict Observation~\ref{obs2}. 

To find this~$t_0$ with the required properties, note that the work between~$\tequal$ and~$t^>$ on both~$i$ and~$j$ is equal (as both are in~$\T$). However, we have~$T_i(\tequal) < T_j(\tequal)$ and~$T_i(t^>) > T_j(t^>)$. In other words, even though the same amount of work has been done, job~$i$ was cooled less. So, at some point between~$\tequal$ and~$t^>$ we must have~$T_i=0$. Define 
\[t_0 = \max\{t \in (\tequal, t^>): T_i =0 \}\]
as the last moment before~$t^>$ that this occurs. In particular this means~$T_i(t)>0$ for all~$ t \in (t_0,t^>)$. So we can use Claim~\ref{claim:workdone} to find that:

\begin{align*}
    \lefttodo_i(t_0) - \lefttodo_i(t^>) &= \frac{T_i(t^>) - T_i(t_0)- \alpha(t^> - t_0)}{\beta - \alpha} \\
    &> \frac{T_j(t^>) - T_j(t_0)- \alpha(t^> - t_0)}{\beta - \alpha} \\
    &\ge \lefttodo_j(t_0) - \lefttodo_j(t^>),
\end{align*}
where we use~$T_i(t^>) > T_j(t^>)$ and~$T_i(t_0) = 0 \le T_j(t_0)$ in the first inequality. Using~$\lefttodo_i(t^>) =
\lefttodo_j(t^>)$ as~$t^> \in \T$, we get what we were after, namely~$\lefttodo_i(t_0) > \lefttodo_j(t_0)$. 
\end{claimproof}

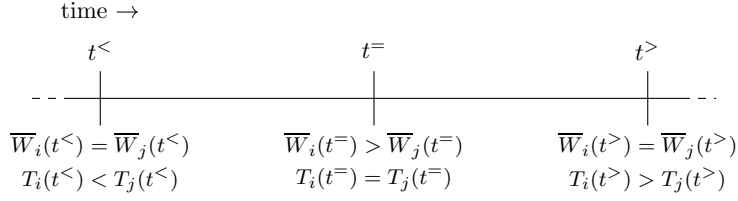
\begin{figure*}[h!]
    \centering
    \begin{tikzpicture} [scale = .9, transform shape,
]

\draw[dashed] (-1,0) -- (-.5,0);
\draw[color = black] (-0.5,0) -- (8.5,0);
\draw[dashed] (8.5,0) -- (9,0);

 \node[] at (0,1.3) 
    {time $\rightarrow$};
 \node[] at (0,.7) {$\tequal$};
 
 \node[] at (0,-.7) {\small $\lefttodo_i(\tequal) = \lefttodo_j(\tequal) $};
\node[] at (0,-1.2) {\small $T_i(\tequal) < T_j(\tequal) $};
\draw[] (0,-.4) -- (0,.4);

\node[] at (8,.7) {$t^>$};
\draw[] (8,-.4) -- (8,.4);
\node[] at (8,-.7) {\small $\lefttodo_i(t^>) = \lefttodo_j(t^>) $};
\node[] at (8,-1.2) {\small $T_i(t^>) > T_j(t^>) $};

\node[] at (4,.7) {$t^=$};
\draw[] (4,-.4) -- (4,.4);
\node[] at (4,-.7) {\small $\lefttodo_i(t^=) > \lefttodo_j(t^=) $};
\node[] at (4,-1.2) {\small $T_i(t^=) = T_j(t^=) $};


\end{tikzpicture}
    \caption{The properties and existence of $t^<$ and $t^>$ follow from Claims~\ref{claim:2} and~\ref{claim:3} (see also Figure~\ref{fig:timeline1}). Claim~\ref{claim:4} ensures the existence of $t^= = \max\{t \in (t^<,t^>): T_i(t) = T_j(t)\}$. The property  $\lefttodo_i(t^=) > \lefttodo_j(t^=)$ follows from Claim~
\ref{claim:5}. }
    \label{fig:timeline2}
\end{figure*}

See Figure~\ref{fig:timeline2} for an illustration. 

\begin{claim}
    We have ~$T'_j(t^=) \le T_j(t^=) = T_i(t^=)$.
\end{claim}
\begin{claimproof}
    Recall everything that we know about~$t'$: job~$j$ overheats directly after~$t'$ in schedule~$S'$,~$\lefttodo_i(t') < \lefttodo_j(t')$ and~$T_j'(t) > 0$ for all~$t \in (\tequal,t')$ (see the proof of Claim~\ref{claim:t<exists}). 
    
First, note that~$t' > t^>$: we would otherwise have~$\lefttodo_i(t^=) > \lefttodo_j(t^=)$ and~$\lefttodo_i(t') < \lefttodo_j(t')$ where~$t' < t^>$, implying (using Observation~\ref{obs:1}) that there exists a~$t \in (\tequal,t^>)$ with~$t \in \T$ contradicting Observation~\ref{obs2}. 

Since we have 
$\lefttodo_i(t^=) > \lefttodo_j(t^=)$ and~$\lefttodo_i(\tequal) = \lefttodo_j(\tequal)$, this means that more work has been done on~$j$ than on~$i$ between~$\tequal$ and~$t^=$ in~$S$. As~$S'_j(t) = S_i(t)$ for all~$t \ge \tequal$, this also means that more work has been done on~$j$ in~$S$ than on~$j$ in~$S'$ between~$\tequal$ and~$t^=$. Combining this with~$T'_j(t)>0$ for all~$t \in [t^<, t']$ and~$t' > t^> > t^=$ we get that (using Claim~\ref{claim:workdone}): 
\begin{align*}
    T'_j(t^=) &= T_j(\tequal) + \int_{\tequal}^{t^=} \alpha(1 - S_i(t)) + \beta S_i(t) \,dt\\
    &\le T_j(\tequal) + \int_{\tequal}^{t^=} \alpha(1 - S_j(t)) + \beta S_j(t) \,dt\\
    &\le T_j(t^=)\\
    &= T_i(t^=).
\end{align*}
\end{claimproof}

We are now ready to prove our contradiction. Since~$T'_j(t^=) \le T_i(t^=)$, it means that, since~$S'_j(t) = S_i(t)$ for all~$t \ge \tequal$, we have~$T'_j(t) \le T_i(t)$ for all~$t > t^=$. In particular, this would imply that job~$i$ would also violate its temperature constraints just after time~$t'$ (as~$j$ does in~$S'$). However, this contradicts that~$S$ is a feasible schedule and finishes the proof of the lemma under the initial assumption \eqref{eq:tempnot1}. We will next reduce the general case to this special case. 

Note that for any linear ordering~$<$ of~$J$, the set of feasible completion time vectors 
\[\mathcal{C}^<:=\left\{ C(S'):\begin{aligned}& S\text{ feasible normal, }\\& i<j\Rightarrow C_i(S')\leq C_j(S')~\forall i,j\in J\end{aligned}\right\}\]
is a polyhedron by Lemma \ref{lem:norm}. As such,~$\mathcal{C}^<$ is a closed set. 

Consider the feasible normal schedule~$S$ that is the subject of the lemma, let~$C:=C(S)\in \R^n$ and let~$C'\in \R^n$ denote~$C$ with~$i$-th and~$j$-th entries switched. To finish the proof, we need to argue that there exists a feasible normal schedule~$S'$ so that~$C(S')=C'$.

Let~$S^\gamma$ be the normal schedule that is determined by
$$S^\gamma(\gamma t)=\frac{S_k(t)}{\gamma}\text{ for all }t\in \R_+.$$
We have~$C(S^\gamma)=\gamma C$. If~$\gamma>1$, then ~$S^\gamma$ is feasible; indeed assumption \eqref{eq:tempnot1} will then hold for~$S^\gamma$, as in such a schedule each job is operated upon in a strictly more relaxed fashion compared to~$S$. Hence, for each~$\gamma>1$, there exists a normal schedule~$(S^\gamma)'$ in which the completion times of~$i,j$ are switched compared to~$S^\gamma$. Then~$C((S^\gamma)')=\gamma C'$. Clearly
$$\lim_{\gamma\downarrow 1} C((S^\gamma)')=\lim_{\gamma\downarrow 1} \gamma C'=C'$$
For any linear order~$<$ of~$J$ so that~$C'\in \mathcal{C}^<$, we have~$\gamma C' = C((S^\gamma)')\in \mathcal{C}^<$. As~$\mathcal{C}^<$ is a closed set, it follows that~$C'\in \mathcal{C}^<$. By definition of~$\mathcal{C}^<$, there exists a feasible normal schedule~$S'$ so that~$C(S')=C'$, as required. \endproof
 In the proof of \Cref{thm:fixedvolgorde}, we have showed how to swap the paths of two jobs when \(C_j < C_i\) while \(p_i < p_j\). This swapping can be done irrespective of the machines that process the jobs - hence, we can extend the theorem to deal with multiple machines. 
 \begin{theorem}
\label{thm:fixedvolgorde_multi}
Given an instance of~$Qm^{sum}$, there is an optimal normal schedule where the jobs are completed in order from smallest to largest processing times. 
\end{theorem}
 Clearly, the LP can be extended to cater for multiple machines \(m\).
Combining this leads to \Cref{cor:multiplemachines}.
\begin{corollary}
\label{cor:multiplemachines}
Any instance of~$Qm^{sum}$ has an optimal solution where the jobs are completed in order from smallest to largest processing times, and this optimal schedule can be found in polynomial time. 
\end{corollary}

\section{Minimizing the Makespan}
\label{sec:makespan}
In this section, we focus on problem~$Qm^{max}$. We will show that this problem can be solved in linear time, even when the heating and cooling rates are job-dependent. Let~$m$ be the number of identical machines, and for every job $j \in J$, we let~$\alpha_j < 0$ be its cooling rate and~$\beta_j > 0$ be its heating rate. Notice that we are still allowed to assume that the maximal temperature is~$1$, else we can normalize the heating and cooling rates by dividing them by the maximal temperature. First define:
\begin{eqnarray*}
    q_j := \text{minimal makespan of a schedule that corre-} \\ \text{sponds to an instance solely containing job } j.
\end{eqnarray*}

\begin{claim}
    \[q_j = \begin{cases}p_j\left(1- \frac{\beta_j}{\alpha_j}\right) + \frac{1}{\alpha_j} &\text{if~$p_j\beta_j \le 1$},\\p_j &\text{otherwise.} \end{cases}\]
\end{claim}

\begin{proof}
    As Theorem~\ref{thm:natural=normal} applies, there is an optimal schedule that is normal. Hence, in this schedule the job is processed at a constant rate~$s_j$ throughout the schedule and we have~$q_js_j = p_j$. There are now two cases, depending on whether or not the job will overheat if we process it at maximum rate~$s_j=1$. If it doesn't (i.e. if $\beta_jp_j \le 1$) then we have~$q_j = p_j$. Otherwise, we observe that in a schedule with minimum makespan, the job has temperature equal to 1 when its processing completes. So we get: 
    \begin{align*}
        && T_j(q_j) &=1 \\
        \Leftrightarrow && \beta_jp_j + \alpha_j (q_j - p_j) &= 1 \\
        \Leftrightarrow && p_j \left(1- \frac{\beta_j}{\alpha_j}\right) + \frac{1}{\alpha_j} &= q_j.
    \end{align*}
\end{proof}

We are now able to solve this problem in linear time. 

\begin{theorem} The minimum makespan can be determined in linear time as 
\[C_{\max} = \max \left\{ \max_{j \in J} \{q_j\}, \frac{\sum_{j \in J}p_j}{m}\right\}.\] 
\end{theorem}
\begin{proof}
    Notice that \[C_{\max} \ge \max \left\{ \max_{j \in J} \{q_j\}, \frac{\sum_{j \in J}p_j}{m}\right\},\] as each of the terms is a lowerbound for the makespan: any instance containing job~$j$ takes at least~$q_j$ time units to complete, and as~$\sum_{j \in J}p_j/m$ represents the total processing time required divided by the number of machine available, this term is also a lowerbound for~$C_{\max}$. 

    Hence, we are left to show that we can find a schedule with the given makespan. Consider the schedule where each job is processed with a rate of~$s_j = p_j /C_{\max}$ during the whole schedule of length~$C_{\max}$. First, it is obvious that each job~$j \in J$ is processed for~$p_j$ time units. Second, we see that for each~$j \in J$, we have~$s_j \le 1$, as~$p_j \le C_{\max}$, so no job is using more than one machine at the same time. 
    Next, we check whether the temperature constraints are satisfied. 
    Indeed, observe that for each job $j \in J$:
   \begin{align*}
       T_j(C_{\max}) &= \alpha_j(C_{\max} - p_j) + \beta_j p_j \\
       &\le \alpha_j(q_j - p_j) + \beta_j p_j\\
       & \le 1,
   \end{align*}
    where the last inequality holds as it is the temperature of $j$ in a feasible normal schedule processing solely $j$ on one machine of makespan $q_j$.
    
    We are left to prove that only~$m$ machines are used at any moment in time. For this we have to check that~$\sum_{j \in J}s_j \le m$. Note that~$C_{\max} \ge \sum_{j \in J}p_j /m$, so 
    \begin{align*}
        \sum_{j \in J}s_j &= \sum_{j \in J}p_j/C_{\max}\\
        &\le \sum_{j \in J}p_j / (\sum_{j \in J}p_j /m) \\
        &= m.
    \end{align*}
\end{proof}

\section{Discussion and open questions}
\label{sec:open}
We have considered a scheduling problem where processing a job raises its temperature, and not processing a job leads to its cooling. The temperature of a job should not exceed a given threshold, and cannot become negative. It is allowed to preempt a job an arbitrary number of times. We are interested in minimizing the makespan, and in minimizing the sum of completion times. Here, we have shown how to solve these problems in polynomial time.

When it comes to open questions: there are many. For the case of minimizing sum of completion times, we list a number of questions and also consider possible generalizations that (seemingly) cannot be handled by the approaches from the previous sections.
\begin{itemize}
\item Is solving the linear program necessary, or does there exist a combinatorial algorithm for minimizing~$\sum_j C_j$?
\item What if jobs have an individual maximum temperature?
\item What if there are weights given for the jobs?
\item What if each job~$j \in N$ has an arbitrary initial temperature~$T_j(0)$?
\item What if jobs have arbitrary heating rates and/or cooling rates?
\end{itemize}

\backmatter

\bmhead{Acknowledgments}
We thank Sem Borst and Mark Peletier for providing pointers to literature and for discussions on the problem. 

The research of Frits Spieksma was partly funded by the NWO Gravitation Project NETWORKS, Grant Number 024.002.003.

C\'eline Swennenhuis has received funding from the European Research Council (ERC) under the European Union’s Horizon 2020 research and innovation programme (grant agreement No 803421, ReduceSearch)


\begin{thebibliography}{99}

\bibitem{Audetetal2004}
Audet, C., J. Brimberg, P. Hansen, S. Le Digabel, N. Mladenovi\'c (2004), {\em Pooling Problem: Alternate Formulations and Solution Methods}, Management Science 50, 761-776.

\bibitem{Baietal2008}
Bai, D., L. Tang, M. Su (2008), {\em A New Machine Scheduling Problem with Temperature Loss}, International Workshop on Knowledge Discovery and Data Mining, pages 662-666, doi  10.1109/WKDD.2008.36

\bibitem{Bampisetal2013}
Bampis, E., D. Letsios, G. Lucarelli, E. Markakis, I. Milis (2013), {\em On multiprocessor temperature-aware scheduling problems}, Journal of Scheduling {\bf 16}, 529-538.

\bibitem{berkim2023}
Bertsimas, D., C.W. Kim (2023), {\em Optimal Control of Multiclass Fluid Queueing Networks: A Machine Learning
Approach}, https://arxiv.org/abs/2307.12405.


\bibitem{BirFun2013}
Birks, M. and S.P.Y. Fung (2013), {\em Temperature aware online algorithms for scheduling equal length jobs}, Theoretical Computer Science {\bf 508}, 54-65.

\bibitem{BirFun2017}
Birks, M. and S.P.Y. Fung (2017), {\em Temperature aware online algorithms for minimizing flow time}, Theoretical Computer Science {\bf 661}, 18-34.

\bibitem{boonetal2011}
Boon, M., R. van der Mei, E. Winands (2011), {\em Applications of polling systems}, Surveys in Operations Research and Management Science 16, 67-82.

\bibitem{borbox2018}
Borst, S., O. Boxma (2018), {\em Polling: past, present, and perspective}, TOP 26, 335-369.

\bibitem{Chrobaketal2011}
Chrobak, M., C. D\"urr, M. Hurand, J. Robert (2011), {\em Algorithms for temperature-aware task scheduling in microprocessor systems}, Sustainable Computing-Informatics \& Systems 1, 241-247.

\bibitem{Gupteetal2017}
Gupte, A., S. Ahmed, S. Dey, M. Seok Cheon (2017), {\em Relaxations and discretizations for the pooling problem}, Journal of Global Optimization 67, 631-669.

\bibitem{Kallrath2002}
Kallrath, J. (2002), {\em Planning and scheduling in the process industry}, OR Spectrum 24, 219-250.

\bibitem{Kilic2011}
Kilic, O. (2011), {\em Planning and scheduling in process industries considering industry-specific characteristics}, PhD thesis of the University of Groningen, available at https://pure.rug.nl/ws/portalfiles/portal/ 2538168/9complete.pdf.

\bibitem{kruskal1968}
Kruskal, J.B. (1968), {\em Work-Scheduling Algorithms: A Nonprobabilistic Queuing Study (with Possible Application to No. 1 ESS)}, The Bell System Technical Journal 48, 2963-2974.

\bibitem{leflamroo2011}
Lefeber, E., S. L\"ammer, J.E. Rooda (2011), {\em Optimal control of a deterministic multiclass queuing system for which several queues can be served simultaneously}, Systems \& Control Letters 60, 524-529.


\bibitem{Marval1998}
Martinelli, F. and P. Valigi (2002), {\em The A scheduling problem for a finite buffer capacity pull-manufacturing system}, in: Proceedings of the 37th IEEE Conference on Decision and Control, pages 2171-2172.

\bibitem{Marval2002}
Martinelli, F. and P. Valigi (2002), {\em The Impact of Finite Buffers on the Optimal Scheduling of
a Single-Machine Two-Part-Type Manufacturing System}, IEEE Transactions on Automatic Control 47, 1705-1710.

\bibitem{Marval2004}
Martinelli, F. and P. Valigi (2004), {\em Dynamic scheduling for a single machine system under different setup and buffer capacity scenarios}, Asian Journal on Control 6, 229-241.

\bibitem{matetal2016}
Matveev, A., V. Feoktistova, K. Bolshakova (2016), {\em On Global Near Optimality of Special Periodic Protocols for Fluid Polling Systems with Setups}, Journal on Optimization Theory and Applications 171, 1055-1070.

\bibitem{Phametal2020}
Pham, T., P. Leyman, and P. de Causmaecker (2020), {\em The intermittent travelling salesman problem}, International Transactions in Operations Research {\bf 27}, 525-548.

\bibitem{Tak2000}
Takagi, H. (2000), {\em Analysis and application of polling models}, in: Performance Evaluation: Origins and Directions, LNCS 1769, Springer, Berlin, pp. 423–442.

\bibitem{Tullemans2019}
Tullemans, R. (2019), {\em Intermittent travelling salesman problem}, Bachelor Thesis, Eindhoven University of Technology.

\bibitem{Vissem2006}
Vishnevskii, V. and O. Semenova (2006), {\em Mathematical methods to study the polling systems}, Automation and Remote Control 67, 173–220.

\end{thebibliography}
\end{document}